\documentclass[12pt]{article}

\usepackage{times} 
\usepackage{color} 
\usepackage{indentfirst} 
\usepackage{titlesec} 
\usepackage{caption2} 
\usepackage{latexsym,bm} 
\usepackage{amscd,amsfonts,amsmath,amssymb,amsthm,amsxtra} 
\usepackage{mathrsfs}
\usepackage{epsfig,graphicx} 
\usepackage[all]{xy} 
\usepackage[dvipdfm, 
            pdfstartview=FitH,
            CJKbookmarks=true,
            bookmarksnumbered=true,
            bookmarksopen=true,
            colorlinks, 
            pdfborder=001, 
            linkcolor=blue,
            anchorcolor=blue,
            citecolor=blue
            ]{hyperref} 
\usepackage[pagewise]{lineno}\linenumbers 
\usepackage{remreset}
\makeatletter\@removefromreset{figure}{section}\makeatother

\renewcommand{\captionlabeldelim}{.~}

\numberwithin{equation}{section}

\newcaptionstyle{1}{%
  \usecaptionmargin\captionfont%
  \onelinecaption%
  {{\bfseries\captionlabelfont\captionlabel\captionlabeldelim}
    \captiontext}%
  {{\centering\bfseries\captionlabelfont\captionlabel\par}%
    \captiontext}}
    \newcaptionstyle{2}{%
  \usecaptionmargin\captionfont%
  {\centering\bfseries\captionlabelfont\captionlabel\par}
   \onelinecaption{\captiontext}{\captiontext}}

\def\cpts{\captionstyle}

\def\htg{\hypertarget} 
\def\hlk{\hyperlink} 

\def\lb{\label}
\def\rf{\ref}

\def\rT{Theorem \rf}

\def\rP{Proposition \rf}
\def\rR{Remark \rf}
\def\rD{Definition \rf}

\def\rF{Fig \rf}

\theoremstyle{plain}
\newtheorem{Thm}{Theorem}[section]
\newtheorem{Cor}[Thm]{Corollary}

\newtheorem{Lem}[Thm]{Lemma}
\newtheorem{Prop}[Thm]{Proposition}
\theoremstyle{definition}
\newtheorem{Def}[Thm]{Definition}
\newtheorem{Ex}[Thm]{Example}
\newtheorem{Not}[Thm]{Notation}
\newtheorem{Rem}[Thm]{Remark}
\theoremstyle{remark}
\newtheorem*{Pf}{Proof}

\def\bT{\begin{Thm}}
\def\eT{\end{Thm}}
\def\bC{\begin{Cor}}
\def\eC{\end{Cor}}
\def\bD{\begin{Def}}
\def\eD{\end{Def}}
\def\bE{\begin{Ex}}
\def\eE{\end{Ex}}
\def\bL{\begin{Lem}}
\def\eL{\end{Lem}}
\def\bN{\begin{Not}}
\def\eN{\end{Not}}
\def\bP{\begin{Prop}}
\def\eP{\end{Prop}}
\def\bR{\begin{Rem}}
\def\eR{\end{Rem}}
\def\bPf{\begin{Pf}}
\def\ePf{\QED\end{Pf}}
\def\ePff{\end{Pf}}

\def\ba{\begin{array}}
\def\ea{\end{array}}
\def\be{\begin{equation}}
\def\ee{\end{equation}}
\def\bea{\begin{eqnarray}}
\def\eea{\end{eqnarray}}
\def\balgd{\begin{aligned}}
\def\ealgd{\end{aligned}}
\def\bfg{\begin{figure}}
\def\efg{\end{figure}}

\def\bcc{\begin{center}}
\def\ecc{\end{center}}

\def\<{\langle}
\def\>{\rangle}
\def\bc{\bigcup}

\def\cd{\cdots}

\def\bs{\setminus}
\def\es{\emptyset}

\def\ev{\equiv}
\def\sbs{\subset}

\def\ex{\exists}
\def\fa{\forall}
\def\gsl{\geqslant}
\def\lsl{\leqslant}
\def\lb{\label}

\def\ts{\times}

\def\br{\breve}

\def\lra{{\Leftrightarrow}}
\def\ra{{\rightarrow}}
\def\Ra{{\Rightarrow}}
\def\lt{\left}
\def\rt{\right}

\def\mt{\mapsto}

\def\aa{{\alpha}}

\def\ga{{\gamma}}

\def\om{{\omega}}

\def\ep{{\epsilon}}
\def\lm{{\lambda}}

\def\dl{{\delta}}
\def\sg{{\sigma}}

\def\Sg{{\Sigma}}

\def\vf{{\varphi}}

\def\N{{\mathbb N}}

\def\R{{\mathbb R}}
\def\bS{{\mathbb S}}
\def\Z{{\mathbb Z}}
\def\cA{{\mathcal A}}
\def\cB{{\mathcal B}}

\def\cN{{\mathcal N}}
\def\cO{{\mathcal O}}
\def\cP{{\mathcal P}}

\def\sP{\mathscr P}

\def\rd{\mathrm d}

\def\hs{\hspace}

\def\nid{\noindent}
\def\nnb{\nonumber}

\def\tr{\textrm}
\def\mr{\mathrm}
\def\tb{\textbf}
\def\ti{\textit}
\def\dsty{\displaystyle}

\def\bmpg{\begin{minipage}}
\def\empg{\end{minipage}}
\def\hlkh{\hlk{H0}{(H0)}-\hlk{H2}{(H2)}}

\def\QEDclosed{\mbox{\rule[0pt]{1.3ex}{1.3ex}}} 
\def\QED{\QEDclosed} 

\def\td#1{\tilde{#1}}

\def\td#1{\tilde{#1}}

\begin{document}
\title{Multiple Rotation Type Solutions
for \\ Hamiltonian Systems on $T^\ell\ts\R^{2n-\ell}$}

\author{Hui Qiao\\
\footnotesize{{\it Chen Institute of Mathematics and LPMC}}\\
\footnotesize{{\it Nankai University, Tianjin 300071}}\\
\footnotesize{{\it The People's Republic of China}}\\
\footnotesize{{\it E-mail:} qiaohuimath@mail.nankai.edu.cn}}
\date{}

\maketitle{}

\begin{abstract}

{\it This paper deals with multiplicity of rotation type solutions
for Hamiltonian systems on $T^\ell\ts\R^{2n-\ell}$. It is proved
that, for every spatially periodic Hamiltonian system, i.e., the
case $\ell=n$, there exist at least $n+1$ geometrically distinct
rotation type solutions with given energy and rotation vector. It is
also proved that, for a class of Hamiltonian systems on
$T^\ell\ts\R^{2n-\ell}$ with $1\lsl\ell\lsl 2n-1$ but $\ell\neq n$,
there exists at least one periodic solution or $n+1$ rotation type
solutions on every contact energy hypersurface.}

\end{abstract}

{\bf 2010 MSC}: 37J45

{\bf Key words}: Hamiltonian system; \hlk{rots}{rotation type
solution}; \hlk{cttt}{contact type}; \hlk{sympldil}{symplectic
dilation}; \hlk{prdext}{periodic extension}

\setcounter{equation}{0}
\section{Introduction and main result}

\nocite{BW97JDE,BW97MZ,kcC89,CZ83,pF92JDE,pF92TMA,cG91,HZ,mJ94,fJ94PLMS,pR88,aW79}

Existence and multiplicity of periodic solutions of spatially
periodic Hamiltonian systems had been studied in
\cite{BW97JDE,kcC89,CZ83,pF92JDE,cG91,fJ94PLMS,pR88} and some
references therein. This paper is motivated by \cite{pF92TMA}. We
make the following assumptions about Hamiltonian function $H$

\htg{H0}{(H0)}\quad $H\in C^2(\R^{2n},\R)$.

\htg{H1}{(H1)}\quad $H(p,q+m)=H(p,q),\ \fa
(p,q)\in\R^{2n-\ell}\ts\R^\ell,\ m\in\Z^\ell,\ 1\lsl\ell\lsl 2n$.

\htg{H2}{(H2)}\quad $\ex\ \mu,r>0$ such that \[0<\mu H(p,q)\lsl
p\cdot H_p(p,q),\ \fa (p,q)\in\R^{2n-\ell}\ts\R^\ell,\ |p|\gsl r.\]

\nid Spatially periodic Hamiltonian functions are the special case
of $\ell=n$ in \hlk{H1}{(H1)}-\hlk{H2}{(H2)}. If $\ell=2n$, then $H$
only satisfies \hlk{H0}{(H0)} and \hlk{H1}{(H1)} and can be
considered as a Hamiltonian function on the torus $T^{2n}$, on which
Hamiltonian systems had been studied in \cite{BW97MZ,kcC89,CZ83}. If
$1\lsl\ell\lsl2n-1$, an alternative result had been obtained by M.
Y. Jiang in \cite{mJ94}.

We look for solutions of \be\lt\{\balgd
  \dot{z}(t)&=JH'(z(t)),\\
  z(T)&=z(0)+(0,k),\\
  H(z(t))&\ev M,\ealgd\rt.\lb{feH}\ee with prescribed energy $M$
and $k\in\Z^\ell$. A pair $(z,T)$ satisfying (\rf{feH}) with
$k\neq0$ is called a \htg{rots}{{\tb{\ti{rotation type solution}}}}
with \htg{rotv}{{\tb{\ti{rotation vector}}}} $k$. Two solutions
$(z_1,T_1)$ and $(z_2,T_2)$ of (\rf{feH}) are called
\htg{gd}{{\tb{\ti {geometrically distinct}}}} if \be\pi\circ
z_1(\R)\neq\pi\circ z_2(\R),\lb{gd}\ee where $\pi$ denotes the
covering projection
$\R^{2n-\ell}\ts\R^\ell\ra\R^{2n-\ell}\ts\R^\ell/\Z^\ell$. Denote by
$\td\cP(M,k)$ the set of geometrically distinct solutions of
(\rf{feH}). For the constant $r>0$ in \hlk{H2}{(H2)}, set \be
M^\ast:=\max\{H(p,q)\ |\ |p|\lsl r,q\in\R^\ell\}.\lb{defM*}\ee The
main result in this paper is

\bT\lb{rotfe} Assume that $H:\R^{2n}\ra\R$ satisfies \hlkh~with
$\ell=n$. For every $M>M^\ast$ and $k\in\Z^n\bs\{0\}$, we have
$^\#\td\cP(M,k)\gsl n+1$.\eT

\nid As an application of Theorem 1.2 of \cite{mJ94}, we have

\bT\lb{alter} Assume that $H:\R^{2n}\ra\R$ satisfies \hlkh~with
$1\lsl\ell\lsl 2n-1$ but $\ell\neq n$. For $M>M^\ast$, if
$H^{-1}(M)$ is of contact type then $^\#\td\cP(M,0)\gsl1$ or
$^\#\td\cP(M,k_0)\gsl n+1$ for some $k_0\in\Z^\ell\bs\{0\}$.\eT

\bR Assumptions of \rT{rotfe} are the same with that of Theorem 0.1
of \cite{pF92TMA} which states that $^\#\td\cP(M,k)\gsl n$, while
\rT{rotfe} gives one more solution of (\rf{feH}). When $n=1$,
Theorem 0.1 of \cite{pF92TMA} gives one rotation type solution,
while \rT{rotfe} gives two. For example, consider the simple
pendulum equation $\ddot q+\sin(q)=0$ with Hamiltonian function
$H(p,q)=\frac{1}{2}p^2-\cos(q)$. For every $M>1$, the energy surface
$H^{-1}(M)$ consists of two disjoint curves $z_1$ and $z_2$, each
one is a rotation type solution. See \rF{fig.1}.\eR

The idea of the proof of \rT{rotfe} is as follows. In
\cite{pF92TMA}, P. Felmer modified $H$ quadratically to obtain a new
Hamiltonian function $\bar H$ which is regular near $\bar
H^{-1}(M)=H^{-1}(M)$, and then obtained solutions $(\bar x,\bar T)$
of \be\lt\{\balgd\dot x(t)+(0,k)&=\bar TJ\bar
H'(x(t)+t(0,k)),\\ x(1)&=x(0),\\
\bar H(x(t)+t(0,k))&\ev M.\ealgd\rt.\lb{febH}\ee A reparametrization
of $\bar x(t)+t(0,k)$ gives a solution of (\rf{feH}). This proves
the existence. {\bf For the multiplicity proof, instead of
Ljusternik-Schnirelmann category argument in \cite{pF92TMA}, we
adopt a different method:} Recall that Rabinowitz' idea to construct
a Hamiltonian function $H_\Sg$ for a starshaped hypersurface $\Sg$
with respect to the origin is
to set \begin{align} H_\Sg(x)&=1,\ \fa x\in\Sg,\nnb\\
H_\Sg(\lm x)&=\vf(\lm),\ \fa x\in\Sg,\ \lm\gsl0,\nnb\end{align}
where $\vf$ is some well chosen function. In this paper, we use a
modification of this idea to define a Hamiltonian function $\hat H$
which is periodic in all components and possesses $H^{-1}(M)$ as a
component of one of its regular energy hypersurfaces. From a
solution of (\rf{febH}), we obtain a solution $(\hat x,\hat T)$ of
\be\lt\{\balgd\dot
x(t)+(0,k)&=\hat TJ\hat H'(x(t)+t(0,k)),\\
x(1)&=x(0).\ealgd\rt.\lb{fphH}\ee For the number $\hat T$, following
the idea of \cite{CZ83}, we obtain $2n+1$ solutions of (\rf{fphH}).
{\bf These $2n+1$ solutions are classified into two categories and
give at least $n+1$ geometrically distinct solutions of (\rf{feH}).}

To prove \rT{alter}, define a Hamiltonian function $\br H$ which
possesses $H^{-1}(M)$ as a regular energy hypersurface, is periodic
in $q$ and equals to its maximum outside an open neighborhood of
$H^{-1}(M)$ in $\R^{2n}$. By Theorem 1.2 of \cite{mJ94}, we obtain a
solution of (\rf{fphH}) with $k$ equals to some $k_0\in\Z^\ell$. If
$k_0\neq0$, extending $\br H$ periodically in $p$ and repeating the
multiplicity proof of \rT{rotfe}, we obtain $n+1$ geometrically
distinct solutions of (\rf{feH}).

This paper is divided into four sections. In Section 2, we prepare
some preliminaries and prove the existence of (\rf{fphH}). In
Section 3, we follow the idea of \cite{CZ83} to obtain the
multiplicity of (\rf{fphH}). Here, special care in the proof of the
Corollary of Theorem 5 of \cite{CZ83} has to be taken. The existence
proof of (\rf{fphH}) in Section 2.3 yields the existence of rest
points. Then, assuming the number of them is finite, they form a
Morse decomposition. Hence, the existence proof of (\rf{fphH}) is
necessary. Finally in Section 4, we prove \rT{rotfe} and \rT{alter}.

In this paper, we denote by $x\cdot y$ the standard inner product of
$x$ and $y$ in Euclidean space, and $S^d$ the unit sphere in
$\R^{d+1}$.

\setcounter{equation}{0}
\section{Preliminaries}

\subsection{Some properties of Hamiltonian functions $H$ satisfying
(H0)-(H2)}

By \hlkh, the energy hypersurface $H^{-1}(M)$ is bounded in $p$ and
periodic in $q$.

\bP\lb{prptH=M1} {\bf(cf. Lemma 1.1 of \cite{pF92TMA})} For every
$M>M^\ast$, the following holds.

(i) There exists two constants $r'=r'(H,M),r''=r''(H,M)>r$ such that
\be r'\lsl|p|\lsl r'',\ \fa (p,q)\in H^{-1}(M).\lb{pbdd}\ee

(ii) Given $p\in\R^{2n-\ell}\bs\{0\}$, there exists a unique
periodic function $\sg=\sg_p\in C^2(\R^\ell,\R_+)$ such that \be
H(\sg(q)p,q)=M,\ \sg(q+m)=\sg(q),\ \fa q\in\R^\ell,\
m\in\Z^\ell.\lb{punq}\ee

(iii) The energy hypersurface $H^{-1}(M)$ is periodic in $q$ in the
sense of \be H^{-1}(M)=H^{-1}(M)+\{0\}\ts\Z^\ell.\lb{qprd}\ee \eP

\bPf Conclusions of this proposition are the same with that of Lemma
1.1 of \cite{pF92TMA} except the upperbound in (\rf{pbdd}) whose
proof is given below. By \hlk{H2}{(H2)}, we have \be H(p,q)\gsl
a|p|^\mu,\ \fa (p,q)\in\R^{2n-\ell}\ts\R^\ell,\ |p|\gsl
r,\lb{grtH}\ee where \be
a:=\min_{|p|=r}\frac{H(p,q)}{|p|^\mu}>0.\lb{defa}\ee Set \be
r':=\min\{|p|\ |\ (p,q)\in H^{-1}(M)\},\
r'':=\lt(\frac{M}{a}\rt)^{\frac{1}{\mu}}.\lb{pbd}\ee By (\rf{defM*})
and (\rf{grtH})-(\rf{pbd}), we obtain (\rf{pbdd}). \ePf

\bR In (H2) of \cite{pF92TMA}, the condition $\mu>1$ is assumed.
However, in the proof of Lemma 1.1 of \cite{pF92TMA}, condition
$\mu>0$ is sufficient for making $f'(\sg)>0$ and $f(\sg)\ra+\infty$
as $\sg\ra+\infty$ in (1.5) of \cite{pF92TMA}.\eR

By (1.9) of \cite{pF92TMA}, there exists a unique $C^2$ function
\be\aa:\R^{2n}_\ast:= \R^{2n}\bs(\{0\}\ts\R^\ell)\ra\R_+\lb{aa},\ee
which is positive homogeneous of degree one in $p$ such that \be
H\lt(\frac{p}{\aa(p,q)},q\rt)\ev M \tr{\quad and\quad
}\aa(p,q+m)=\aa(p,q),\ \fa (p,q)\in\R^{2n}_\ast,\
m\in\Z^\ell.\lb{defaa}\ee In fact, for any $\lm>0$, we have
\[M=H\lt(\frac{\lm p}{\aa(\lm p,q)},q\rt)=H\lt(\frac{p}{\aa(p,q)},q\rt).\]
Then, by the uniqueness of $\aa$, we have \be\frac{\lm}{\aa(\lm
p,q)}=\frac{1}{\aa(p,q)}\ \Ra\ \aa(\lm p,q)=\lm\aa(p,q).\lb{aalm}\ee

\bP\lb{prptH=M2} For every $M>M^\ast$, the energy hypersurface
$H^{-1}(M)$ is $C^2$ diffeomorphic to $S^{2n-\ell-1}\ts\R^\ell$.\eP

\bPf By (\rf{defaa}), we have \be H^{-1}(M)=\aa^{-1}(1)
=\bigcup_{q\in\R^\ell}\{p\in\R^{2n-\ell}\ |\
\aa(p,q)=1\}.\lb{HM=aa1}\ee Define two maps
\begin{align}&F_1:H^{-1}(M)\ra S^{2n-\ell-1}\ts\R^\ell,\
F_1(p,q):=\lt(\frac{p}{|p|},q\rt),\lb{defF1}\\
&F_2:S^{2n-\ell-1}\ts\R^\ell\ra H^{-1}(M),\
F_2(p,q):=\lt(\frac{p}{\aa(p,q)},q\rt).\lb{defF2}\end{align} For any
$(p,q)\in H^{-1}(M)$, we have $\aa(p,q)=1$ and \[F_2\circ
F_1(p,q)=F_2\lt(\frac{p}{|p|},q\rt)
=\lt(\frac{\frac{p}{|p|}}{\aa\lt(\frac{p}{|p|},q\rt)},q\rt)
=\lt(\frac{p}{\aa(p,q)},q\rt)=(p,q).\] This proves $F_2\circ
F_1=id_{H^{-1}(M)}$. Similarly, we have $F_1\circ
F_2=id_{S^{2n-\ell-1}\ts\R^\ell}$. By (i) of \rP{prptH=M1}, we have
$H^{-1}(M)\sbs\R^{2n}_\ast$. Since $\aa\in C^2(\R^{2n}_\ast,\R)$,
then $F_2$ defined by (\rf{defF2}) is a $C^2$-diffeomorphism. \ePf

\bE Consider the Hamiltonian function $H:\R^4\ra\R$ defined by
$H(p_1,p_2,q_1,q_2):=\frac{1}{2}(|p_1|^2+|p_2|^2)-
\cos(q_1)-\cos(q_2)$. For every $M>2$, the energy hypersurface
$H^{-1}(M)$ consists of circles with radius varying periodically
from $\sqrt{2(M-2)}$ to $\sqrt{2(M+2)}$.\eE

\subsection{Symplectic dilation}

In this subsection, we prepare some results for the multiplicity
proof in Section 4.

\bD\lb{wdltn} {\bf(\cite{aW79})} A
\htg{sympldil}{{\tb{\ti{symplectic dilation}}}} of a hypersurface
$S$ in a symplectic manifold $(M,\om)$ is a vector field $\xi:U\ra
TM$ defined on an open neighborhood $U$ of $S$ in $M$ which is
transverse to $S$ and satisfies \be
L_\xi\omega=\omega.\lb{wdlt}\ee\eD

A. Weinsten \cite{aW79} gives the following equivalent
characterization of contact type.

\bP\lb{contact} A hypersurface $S$ in $(M,\omega)$ is of
\htg{cttt}{{\tb{contact type}}} if and only if there exists a
symplectic dilation of $S$.\eP

Let $\xi:U\ra TM$ be a symplectic dilation of a hypersurface $S$ in
$(M,\om)$. For each $x\in U$, there exist $\dl(x)>0$, a neighborhood
$U(x)$ of $x$ in $U$ and a map $\Phi:(-\dl(x),\dl(x))\ts U(x)\ra
U(x)$ such that \be\lt\{\balgd\frac{\rd}{\rd s}\Phi(s,y)&=
\xi(\Phi(s,y)),\ s\in(-\dl(x),\dl(x)),\\
\Phi(0,y)&=y,\ y\in U(x).\ealgd\rt.\lb{defPhi}\ee If there exists a
constant $\dl:=\dl(S)>0$ such that $\dl(x)\ev\dl,\ \fa x\in S$, then
\be\Phi:(-\dl,\dl)\ts S\ra\Phi((-\dl,\dl)\ts
S):=\cN_\dl\lb{prptPhi}\ee is a diffeomorphism. If $\xi$ is $C^N$,
where $N\in\N\cup\{+\infty\}$, then $\Phi$ is $C^N$.

\bD\lb{dlt} We call $\Phi$ a \htg{sympldldil}{{\tb{\ti{symplectic
$\dl$-dilation}}}} of $S$. The $1$-parameter diffeomorphism group
\[\{\phi^s:=\Phi(s,\cdot)\}_{|s|<\dl}\] is also called a symplectic
$\dl$-dilation of $S$.\eD

Let $\{\phi^s\}_{|s|<\dl}$ be a symplectic $\dl$-dilation of a
hypersurface $S$ in $(M,\om)$. A subsequent property of (\rf{wdlt})
is \be(\phi^s)^\ast\om=e^s\om,\ \fa s\in(-\dl,\dl),\lb{phisom}\ee
which is crucial in the proof of (ii) of \rP{keyprop} below. We
refer to page 122 of \cite{HZ} for its proof. The following map \be
\phi^s:\cN_{\frac{1}{3}\dl}\ra\phi^s(\cN_{\frac{1}{3}\dl}),\ \fa
s\in(-\frac{2}{3}\dl,\frac{2}{3}\dl),\lb{phis}\ee is a
diffeomorphism. By (\rf{prptPhi}) and Definition \rf{dlt}, we have
\be\cN_\dl=\bc_{|s|<\dl}S_s,\tr{ where
}S_s:=\phi^s(S).\lb{NdltaSs}\ee Choose $\hat H\in C^2(M,\R)$ such
that \be\hat H(x):=f(s),\ \fa s\in(-\frac{2}{3}\dl,\frac{2}{3}\dl),\
x\in S_s,\lb{defhH}\ee where $f:\R\ra\R$ is a $C^2$ auxiliary
function and satisfies \be\lt\{\balgd f(s)\ev \max_{s\in\R}
f&,&&s\gsl\frac{1}{3}\dl,\\f'(s)>0&,&&|s|<\frac{1}{3}\dl.\ealgd\rt.\lb{prptf}\ee
See \rF{fig.f}. By (\rf{prptPhi}), the smoothness of $\hat H$
follows from that of $f$.

\bP\lb{keyprop} Let $S$ be a contact hypersurface in a symplectic
manifold $(M,\om)$, $\xi:U\ra TM$ be a symplectic dilation of $S$
which generates $\{\phi^s\}_{|s|<\dl}$, $\hat H$ be defined by
(\rf{defhH}) and $X_{\hat H}$ be the Hamiltonian vector field of
$\hat H$. The following holds

(i) For each $s\in(-\frac{1}{3}\dl,\frac{1}{3}\dl)$ and $x\in S_s$,
we have \be\hat H'(x)\cdot\xi(x)=f'(s)>0\lb{hHtrsvs}\ee and $S_s$ is
a regular energy hypersurface of $\hat H$.

(ii) For $s_1, s_2\in(-\frac{1}{3}\dl,\frac{1}{3}\dl)$, if $x_1(t)$
lies on $D_{s_1}$ and solves $\dot x_1(t)=X_{\hat H}(x_1(t))$, then
\[x_2(t):=\phi^{s_2-s_1}(x_1(t))\] lies on $S_{s_2}$
and solves \be\dot x_2(t)= e^{s_2-s_1}\frac{f'(s_1)}{f'(s_2)}X_{\hat
H}(x_2(t)).\lb{dtx2}\ee

(iii) If $F$ is another function possessing $S$ as a regular energy
surface, then\be X_F(x)=(F'(x)\cdot\xi(x))X_{\hat H}(x),\ \fa x\in
S.\lb{XFx}\ee \eP

\bPf For (i), the conclusion follows from differentiating
(\rf{defhH}) with respect to $s$.

For (ii), since $\{\phi^s\}$ is generated by $\xi$, $\hat H$ and
$\hat H\circ\phi^{s_1-s_2}$ are constant on $S_{s_2}$ and regular
near $S_{s_2}$, by (\rf{hHtrsvs}), we have \be(\hat H\circ
\phi^{s_1-s_2})'(x)=\frac{f'(s_1)}{f'(s_2)}\hat H'(x),\ \fa
s_1,s_2\in(-\frac{1}{3}\dl,\frac{1}{3}\dl),\ x\in
S_{s_2}.\lb{hHphis1s2'}\ee Note that \be d\hat H(x_1(t))=d(\hat
H\circ\phi^{s_1-s_2})(x_2(t))d(\phi^{s_2-s_1})(x_1(t)).
\lb{dhHx1}\ee Then, for any vector field $\zeta$ on
$\cN_{\frac{1}{3}\dl}$, we have
\begin{align}&\om\lt(\frac{\rd}{\rd
t}x_2(t),d(\phi^{s_2-s_1})(x_1(t))\zeta(x_1(t))\rt)\nnb\\
=\,&\om\lt(d(\phi^{s_2-s_1})(x_1(t))\frac{\rd}{\rd t}x_1(t)
,d(\phi^{s_2-s_1})(x_1(t))\zeta(x_1(t))\rt)\nnb\\
=\,&e^{s_2-s_1}\om\lt(\frac{\rd}{\rd t}x_1(t),\zeta(x_1(t))\rt)
=-e^{s_2-s_1}d\hat H(x_1(t))\zeta(x_1(t))\nnb\\
=\,&-e^{s_2-s_1}d(\hat H\circ\phi^{s_1-s_2})(x_2(t))
d(\phi^{s_2-s_1})(x_1(t))\zeta(x_1(t))\nnb\\
=\,&-e^{s_2-s_1}\frac{f'(s_1)}{f'(s_2)}d\hat H(x_2(t))
d(\phi^{s_2-s_1})(x_1(t))\zeta(x_1(t))\nnb\\
=\,&\om\lt(e^{s_2-s_1}\frac{f'(s_1)}{f'(s_2)}X_{\hat
H}(x_2(t)),d(\phi^{s_2-s_1})(x_1(t))\zeta(x_1(t))\rt).\lb{omdtx2}\end{align}
The second equality follows from (\rf{phisom}). The fourth equality
follows from (\rf{dhHx1}). The fifth equality follows from
(\rf{hHphis1s2'}). Since $\om$ is non-degenerate and
$d(\phi^{s_2-s_1})(x_1(t)): T_{x_1(t)}M\ra T_{x_2(t)}M$ is
isomorphic, by (\rf{omdtx2}), we have
\[\frac{\rd}{\rd t}x_2(t)=e^{s_2-s_1}\frac{f'(s_1)}{f'(s_2)}
X_{\hat H}(x_2(t)).\]

For (iii), since $F$ and $\hat H$ are constant on $S$ and regular
near $S$, then the conclusion follows from (\rf{hHtrsvs}) with
$s=0$. \ePf

By elementary calculations, we have

\bP\lb{lft} Let $\td M$ be a covering space of $M$ with covering
projection $\pi:\td M\ra M$. The following holds.

(i) If $\om$ is a symplectic form on $M$, then $\td\om$ defined by
\be\td\om:=\pi^\ast\om\lb{deftdom}\ee is a symplectic form on $\td
M$.

(ii) If $\xi:U\sbs M\ra TM$ is a symplectic dilation of a
hypersurface $S$ in $(M,\om)$, then $\td\xi$ defined by
\be\td\xi(\td x):=d\pi(\td x)^{-1}\xi(\pi(\td x)),\ \fa\td x\in\td
U:=\pi^{-1}(U),\lb{deftdxi}\ee is a symplectic dilation of $\td
S:=\pi^{-1}(S)$.

(iii) If $\{\phi^s\}_{|s|<\dl}$ is a symplectic $\dl$-dilation of a
hypersurface $S$ in $(M,\om)$, there exists a symplectic
$\dl$-dilation $\{\td\phi^s\}_{|s|<\dl}$ of $\td S$ such that \be
\pi\circ\td\phi^{s}(\td x)=\phi^{s}\circ\pi(\td x),\ \fa
s\in(-\dl,\dl),\ \td x\in\td S.\lb{lftphi}\ee \eP

\subsection{Hamiltonian function $\hat H$ for the multiplicity proof}

In this subsection, we list some properties of $\bar H$ given by P.
Felmer \cite{pF92TMA} for the fixed energy problem (\rf{febH}). For
the fixed periodic problem (\rf{fphH}), we define a new Hamiltonian
function $\hat H$ similar to (\rf{defhH}), which possesses
$H^{-1}(M)$ as a component of its regular energy hypersurface.

\bP\lb{prptbH} {\bf(cf. (1.9)-(1.16) of \cite{pF92TMA})} For every
$M>M^\ast$, there exists a function $\bar H:\R^{2n}\ra\R$ such that

(i) $\bar H$ satisfies \hlk{H0}{(H0)}, \hlk{H1}{(H1)} and
\[\mr{(\bar H2)}\hs{0.3cm}0<2\bar H(x)\lsl p\cdot\bar H_p(x),\ \fa
x=(p,q)\in\R^{2n}_\ast.\]

(ii) There exists constants $a_1,a_2,a_3,a_4>0$ such that
\begin{align} a_1|p|^2\lsl\,&\bar H(x)\lsl a_2|p|^2,\
\fa x=(p,q)\in\R^{2n-\ell}\ts\R^\ell,\lb{bH}\\
2a_1|p|\lsl\,&|\bar H'(x)|\lsl a_3(|p|+1),\fa
x=(p,q)\in\R^{2n-\ell}\ts\R^\ell,\lb{bH'}\\
\,&|\bar H''(x)|\lsl a_4,\ \fa
x=(p,q)\in\R^{2n}_\ast.\lb{bH''}\end{align}

(iii) $\bar H$ is regular near $\bar H^{-1}(M)=H^{-1}(M)$ and \be
p\cdot\bar H_p(p,q)=2M>0,\ \fa (p,q)\in\bar
H^{-1}(M).\lb{eq:2bH=}\ee\eP

\bP\lb{exdlt} Assume that $H:\R^{2n}\ra\R$ satisfies \hlkh~with
$\ell=n$. For $M^\ast$ defined by (\rf{defM*}) and every $M>M^\ast$,
the energy surface $H^{-1}(M)$ is of contact type in the sense of
\cite{aW79} (see \rD{wdltn} and \rP{contact}).\eP

\bPf Define a vector field \be\xi:\R^{2n}_\ast\ra\R^{2n}_\ast,\
(p,q)\mt(p,0),\ \fa (p,q)\in\R^{2n}_\ast.\lb{xi1}\ee By (i) of
Proposition \rf{prptbH}, we have \[\bar
H'(e^sp,q)\cdot\xi(e^sp,q)=\bar H_p(e^sp,q)\cdot e^sp \gsl2\bar
H(e^sp,q)>0, \ \fa s\in\R,\ (p,q)\in\R^{2n}_\ast.\] This proves that
$\xi$ is a symplectic dilation of $\bar H^{-1}(M)=H^{-1}(M)$. \ePf

\bR The vector field $\xi$ defined by (\rf{xi1}) generates the flow
\be\Phi:\R\ts\R^{2n}_\ast\ra\R^{2n}_\ast,\ (s,p,q)\mt(e^sp,q),\ \fa
s\in\R,\ (p,q)\in\R^{2n}_\ast,\lb{Phi1}\ee which induces
diffeomorphisms \be\Phi:\R\ts\bar
H^{-1}(M)\ra\R^{2n}_\ast\lb{Phi1'}\ee and
\be\phi^s:\R^{2n}_\ast\ra\R^{2n}_\ast,\
(p,q)\mt(e^sp,q).\lb{phis-1}\ee \eR

For $H$ satisfying \hlkh, assume that $H^{-1}(M)$ is contact. By
\rP{prptH=M1}, the hypersurface $H^{-1}(M)/\Z^\ell\sbs
\R^{2n-\ell}\ts\R^\ell/\Z^\ell$ is compact. By Propositions
\rf{lft}-\rf{exdlt}, there exist a constant $\dl>0$ sufficiently
small, two symplectic $\dl$-dilations $\{\phi^s\}$ and
$\{\hat\phi^s\}$ of $H^{-1}(M)$ and $H^{-1}(M)/\Z^\ell$ respectively
such that (\rf{lftphi}) holds. Set
\be\cN_\dl:=\bigcup_{|s|<\dl}D_s,\ {D_s}:=\phi^s(D_0),\
{D_0}:=H^{-1}(M).\lb{NdltaD}\ee

\bP For every $s\in(-\dl,\dl)$, we have \be
D_s=D_s+\{0\}\ts\Z^\ell\lb{Dsprd}\ee and
\be\phi^{-s}(x+(0,m))=\phi^{-s}(x)+(0,m),\ x\in D_s,\
m\in\Z^\ell.\lb{phi-sx0m}\ee\eP

\bPf For any $s\in(-\dl,\dl),\ y\in D_0,\ m\in\Z^\ell$, by
(\rf{qprd}) and (\rf{lftphi}), we have $y+(0,m)\in D_0$ and
\[\pi\circ\phi^s(y+(0,m))=\hat\phi^s\circ\pi(y+(0,m))
=\hat\phi^s\circ\pi(y)=\pi\circ\phi^s(y).\] Since
$\phi^s(y+(0,m))-\phi^s(y)$ is continuous in $s$, there exists
$m'\in\Z^\ell$ independent on $s$ such that
\[(0,m')=\phi^s(y+(0,m))-\phi^s(y)
=\phi^0(y+(0,m))-\phi^0(y)=y+(0,m)-y=(0,m).\] Then $m'=m$ and
\[\phi^s(y+(0,m))=\phi^s(y)+(0,m),\ \fa s\in(-\dl,\dl),\ y\in D_0,\
m\in\Z^\ell.\] This proves (\rf{Dsprd}) and (\rf{phi-sx0m}). \ePf

Define auxiliary functions $f,g:\R\ra\R$ by

\begin{align}
f(s)&:=\lt\{\ba{ccc} 0&,&s\lsl-\frac{1}{3}\dl,\\
              \dsty{\frac{1}{(\frac{1}{3}\dl)^2}(s+\frac{1}{3}\dl)^3
              +\frac{-\frac{1}{2}}{(\frac{1}{3}\dl)^3}(s+\frac{1}{3}\dl)^4}
&,&-\frac{1}{3}\dl\lsl s\lsl0 , \\
              \dsty{\frac{1}{3}\dl+\frac{1}{(\frac{1}{3}\dl)^2}(s-\frac{1}{3}\dl)^3
              +\frac{\frac{1}{2}}{(\frac{1}{3}\dl)^3}(s-\frac{1}{3}\dl)^4}
&,&0\lsl s\lsl\frac{1}{3}\dl, \\
              \dsty{\frac{1}{3}}\dl&,&
              s\gsl\frac{1}{3}\dl,\ea\rt.\lb{deff}\\
g(s)&:=\frac{f'(s)}{e^s}.\lb{defg}\end{align}
 By direct calculation, we have

\bP\lb{auxfg}\quad

(i) $f$ is $C^2$ and satisfies (\rf{prptf}).

(ii) $g$ is $C^1$ and satisfies \be
g(s)\lt\{\balgd=0&,&&|s|\gsl\frac{1}{3}\dl,\\
>0&,&&|s|<\frac{1}{3}\dl,\ealgd\rt.\tr{and
}\lt\{\balgd&\tr{increases
strictly in }s\in[-\frac{1}{3}\dl,-\frac{1}{3}\dl+\dl^+],\\
&\tr{decreases strictly in }s\in[-\frac{1}{3}\dl+\dl^+,
\frac{1}{3}\dl],\ealgd\rt.\lb{prptg}\ee where \be
\dl^+:=\frac{3+\frac{1}{2}\dl-\sqrt{(3+\frac{1}{2}\dl)^2-4\dl}}{2}.
\lb{defdlta+}\ee See \rF{fig.g}. \eP

\bP\lb{pext} {\bf (\htg{prdext}{Periodic extension})} For the vector
$k\in\Z^\ell\bs\{0\}$ in \rT{rotfe}, there exist an $N_0\in\N$
sufficiently large and a function $\hat H\in C^2(\R^{2n},\R)$ such
that

(i) $\hat H$ is periodic in all components in the sense of \be\hat
H(x+(m',m))=\hat H(x),\ \fa x\in\R^{2n},(m',m)\in
N_0\Z^{2n-\ell}\ts\Z^\ell,\lb{pexthH}\ee and possesses bounded
gradient and Hessian, i.e., \be 0<M_1:=\max_{x\in\R^{2n}}|\hat
H'(x)|,\ M_2:=\max_{x\in\R^{2n}}|\hat
H''(x)|<+\infty.\lb{hH'hH''}\ee

(ii) For each $b\in(0,\frac{1}{3}\dl)$, there exists some
$\dl_b\in(-\frac{1}{3}\dl,\frac{1}{3}\dl)$ such that \be\hat
H^{-1}(b)=D_{\dl_b}+N_0\Z^{2n-\ell}\ts\{0\},\lb{hH-1b}\ee where
$D_s$ is defined by (\rf{NdltaD}). Moreover, \be\hat
H'(x)\neq0\tr{\quad iff\quad}0<\hat
H(x)<\frac{1}{3}\dl.\lb{hH'neq0}\ee

(iii) If $(\bar x,\bar T)$ solves (\rf{febH}) with $\bar T\neq0$,
then \be \hat x(s):=\bar x(t(s))+t(s)(0,k)-s(0,k)\lb{defys}\ee
solves (\rf{fphH}) with $\hat T$ defined by \be\hat T:=\bar
T\int^1_0\bar H'(\bar x(t)+t(0,k))\cdot\xi(\bar x(t)+t(0,k))\rd
t\neq0,\lb{defhT}\ee where $\xi$ denotes the symplectic dilation of
$D_0$ which generates $\phi^s$ in (\rf{NdltaD}) and $t(s)$ is
defined by \be\lt\{\balgd\frac{\rd t}{\rd s}&=\hat
T(\bar T\bar H'(\bar x(t)+t(0,k))\cdot\xi(\bar x(t)+t(0,k)))^{-1}\neq0,\\
t(0)&=0,\ealgd\rt.\lb{defts}\ee

(iv) If $\hat x$ solves (\rf{fphH}), then $\hat H(\hat
x(t)+t(0,k))\ev b\in(0,\frac{1}{3}\dl)$. Moreover,

(iv.1) Each element of \be[\hat x]:=\hat
x+N_0\Z^{2n-\ell}\ts\Z^\ell\lb{def[x]}\ee solves (\rf{fphH}).

(iv.2) There exist $\dl_b\in(-\frac{1}{3}\dl,\frac{1}{3}\dl)$ and
$x\in[\hat x]$ solving (\rf{fphH}) such that $x(t)+t(0,k)$ lies on
$D_{\dl_b}$. Then \be w(t):=\phi^{-\dl_b}(x(t)+t(0,k))\in D_0,\ \fa
t\in\R.\lb{defwt}\ee

(iv.3) $\pi\circ w$ is a non-constant $1$-periodic solution of
\be\lt\{\balgd\frac{\rd}{\rd t}\pi\circ
w(t)&=g(\dl_b)\hat TX_{\hat H_\pi}(\pi\circ w(t)),\\
\pi\circ w(1)&=\pi\circ w(0),\ealgd\rt.\lb{hspiw}\ee where $\pi$
denotes the covering projection
$\R^{2n}\ra\R^{2n-\ell}\ts\R^{\ell}/\Z^\ell$, $\hat
H_\pi:\R^{2n-\ell}\ts\R^{\ell}/\Z^\ell\ra\R$ is a Hamiltonian
function defined by \be\hat H_\pi\circ\pi:=\hat H,\lb{hH_pi}\ee
$X_{\hat H_\pi}$ denotes the Hamiltonian vector field of $\hat
H_\pi$ on $(\R^{2n-\ell}\ts\R^{\ell}/\Z^\ell,\om_\pi)$ with
$\om_\pi$ determined by \be\om=\pi^\ast\om_\pi,\lb{om_pi}\ee and
$\om$ denotes the standard symplectic form on $\R^{2n}$.

(iv.4) Denote by $\td\sP$ the set of periods of $\pi\circ w$. There
exists $n(x)\in\N$ such that
\be\td\sP=\frac{1}{n(x)}\Z.\lb{disprd}\ee

(iv.5) \be z(s):=w\lt(t\lt(\frac{s}{T}\rt)\rt)\lb{defzs}\ee solves
(\rf{feH}) with $T$ defined by \be T:=g(\dl_b)\hat
T\int^1_0\frac{\rd t}{H'(w(t))\cdot\xi(w(t))}\neq0,\lb{defT}\ee
where $\hat s:=\frac{s}{T}$ and $t(\hat s)$ is defined by
\be\lt\{\balgd\frac{\rd t}{\rd\hat s}&=\frac{T
H'(w(t))\cdot\xi(w(t))}{g(\dl_b)\hat T}\neq0,\\
t(0)&=0.\ealgd\rt.\lb{defths}\ee\eP

\bPf We carry out the proof in three steps.

{\bf Step 1.} Definition of $\hat H$.

Define a map $\Phi_1$ such that the following diagram commutes,
i.e., \be\Phi_1:=\Phi\circ(id\ts F_2)\circ(F_3\ts
id),\lb{defPhi_1}\ee \be\xymatrix{
(\tr{int}B^{2n-\ell}_{e^\dl}(0)\bs
B^{2n-\ell}_{e^{-\dl}}(0))\ts\R^\ell \ar[r] \ar[d]^{F_3\ts id}
\ar[r]^{\hs{2cm}\Phi_1} & \cN_\dl \ar[r]^{\br H} & \R \\
(-\dl,\dl)\ts S^{2n-\ell-1}\ts\R^\ell \ar[r]^{\hs{0.3cm}id\ts F_2} &
(-\dl,\dl)\ts H^{-1}(M) \ar[u]^{\Phi} & }\lb{cdgr}\ee
\[\xymatrix{ (p,q)
\ar@{|->}[d] \ar@{|->}[r]
& \Phi(\ln|p|,F_2(\frac{p}{|p|},q)) \ar@{|->}[r] & f(\ln|p|)\\
(\ln|p|,\frac{p}{|p|},q) \ar@{|->}[r] &(\ln|p|,\frac{p}{\aa(p,q)},q)
\ar@{|->}[u] & }\] where $\Phi$ is the diffeomorphism (\rf{prptPhi})
with $S=H^{-1}(M)$, $F_2:S^{2n-\ell-1}\ts\R^\ell\ra H^{-1}(M)$ is
the $C^2$-diffeomorphism defined by (\rf{defF2}),
$F_3:\tr{int}B^{2n-\ell}_{e^\dl}(0)\bs B^{2n-\ell}_{e^{-\dl}}(0)\ra
(-\dl,\dl)\ts S^{2n-\ell-1}$ is defined by \be
F_3(p):=(\ln|p|,\frac{p}{|p|}),\ \fa p\in B^{2n-\ell}_{e^\dl}(0)
\bs\bar B^{2n-\ell}_{e^{-\dl}}(0),\lb{defF3}\ee and \be
B^{2n-\ell}_{\rho}(0):=\{p\in\R^{2n-\ell}\ |\ |p|\lsl\rho\},\
\tr{int}B^{2n-\ell}_{\rho}(0):=\{p\in\R^{2n-\ell}\ |\ |p|<\rho\},\
\rho>0.\lb{defB}\ee Then $\Phi_1$ is a diffeomorphism. The
$C^2$-smoothness of $\Phi_1$ follows from that of $F_3,\ F_2$ and
$\Phi$. To continue our study, we define an auxiliary function $\br
H:\R^{2n}\ra\R$ in two cases according to the value of $\ell$ in
\hlk{H1}{(H1)}.

{\bf Case 1.} $1\lsl\ell\lsl2n-2$. By (\rf{defPhi_1}) and
(\rf{cdgr}) with $\dl$ replaced by any $\dl'\in(0,\dl)$, the set
$\R^{2n}\bs\cN_{\dl'}$ has two components $V_{\dl'}$ and $W_{\dl'}$
such that $V_{\dl'}/\Z^\ell$ is compact and $\{0\}\ts\R^\ell$ is
contained in the interior of $V_{\dl'}$. Choose $f$ defined by
(\rf{deff}) and define a new Hamiltonian function $\br H$ by \be\br
H(x):=\lt\{\balgd
  0&,&&x\in V_{\frac{1}{3}\dl},\\
  f(s)&,&&x\in D_s,\ |s|<\frac{2}{3}\dl,\\
  \frac{1}{3}\dl&,&&x\in W_{\frac{1}{3}\dl},\ealgd\rt.\lb{defbrH1}\ee

{\bf Case 2.} $\ell=2n-1$. By (2.65) and (2.66) with $\dl$ replaced
by any $\dl'\in(0,\dl)$, the set $\R^{2n}\bs\cN_{\dl'}$ has three
components $V_{\dl'}$, $W^\pm_{\dl'}$ such that $V_{\dl'}/\Z^\ell$
is compact and $\{0\}\ts\R^\ell$ is contained in the interior of
$V_{\dl'}$. We can also define $\br H$ similar to (\rf{defbrH1}) by
\be\br H(x):=\lt\{\balgd
  0&,&&x\in V_{\frac{1}{3}\dl},\\
  f(s)&,&&x\in D_s,\ |s|<\frac{2}{3}\dl,\\
  \frac{1}{3}\dl&,&&x\in W^-_{\frac{1}{3}\dl}\cup
  W^+_{\frac{1}{3}\dl}.\ealgd\rt.\lb{defbrH2}\ee

It follows from (\rf{pbdd}) that \be D_0=H^{-1}(M)\sbs
B^{2n-\ell}_{r''}(0)\ts\R^\ell,\lb{D0ctdi}\ee Since $\dl>0$ is
sufficiently small, by (\rf{pbdd}), (\rf{NdltaD}) and (\rf{D0ctdi}),
there exists \be r_\dl\in\N\cap(r'',\infty)\lb{rdlta}\ee such that
\be\cN_\dl\sbs
B^{2n-\ell}_{r_\dl}(0)\ts\R^\ell\sbs\tt{int}\cO_\dl\ts\R^\ell,\lb{Ndltactdi}
\ee where the closed set $\cO_\dl$ is defined by \be
\cO_\dl:=\{p=(p_1,p_2,\cd,p_{2n-\ell})\in\R^{2n-\ell}\ |\ |p_i|\lsl
r_\dl+1,\ 1\lsl i\lsl 2n-\ell\}\lb{defOdlta}\ee and
$\tr{int}\cO_\dl$ denotes the interior of $\cO_\dl$. Then $\br H$
defined by (\rf{defbrH1}) or (\rf{defbrH2}) equals to its maximum
$\frac{1}{3}\dl$ on $\R^{2n}\bs(\cO_\dl\ts\R^\ell)$. Choose \be
N_0:=2(r_\dl+1).\lb{defN0}\ee Then we have
\be\R^{2n-\ell}=\bigcup_{m'\in\Z^{2n-\ell}}(\cO_\dl+N_0m')
:=\cO_\dl+N_0\Z^{2n-\ell}\lb{Odlta+N0}\ee and
\be\balgd(\tr{int}\cO_{\dl_b}\ts\R^\ell+(m'_1,0))\cap(\tr{int}\cO_{\dl_b}
\ts\R^\ell+(m'_2,0))=\es,\\ \fa m'_1,m'_2\in
N_0\Z^{2n-\ell},m'_1\neq m'_2.\ealgd\lb{ept}\ee Using $\br H$, we
can define a new Hamiltonian function $\hat H:\R^{2n}\ra\R$ by
\be\hat H(p+N_0m',q)=\br H(p,q),\ \fa p\in\cO_\dl,\ q\in\R^\ell,\
m'\in\Z^{2n-\ell}.\lb{defhH*}\ee It follows from (\rf{Odlta+N0}) and
(\rf{ept}) that $\hat H$ is well-defined. Then $\hat H$ satisfies
(i)-(ii). Here, we define
\[\br F:(\tr{int}B^{2n-\ell}_{e^\dl}(0)\bs
B^{2n-\ell}_{e^{-\dl}}(0))\ts\R^\ell\ra\R,\ (p,q)\mt f(\ln|p|).\]
Then in the diagram (\rf{cdgr}) we have $\br F=\br H\circ\Phi_1$.
Since $f$ and $\Phi_1$ are $C^2$, then $\br F$ and $\br H$ are
$C^2$. Hence $\hat H$ is $C^2$. The periodicity of $\hat H$ in $q$
follows from (\rf{Dsprd}). The most new important property of $\hat
H$ is its periodicity in $p$.

{\bf Step 2.} Proof of (iii).

Since $\xi$ is a symplectic dilation of $D_0$, $\bar H$ is regular
near $\bar H^{-1}(M)=H^{-1}(M)=D_0$ and $\bar H(\bar x(t)+t(0,k))\ev
M$, then $\xi$ is transverse to $D_0$ and \be\bar H'(\bar
x(t)+t(0,k))\cdot\xi(\bar x(t)+t(0,k))\neq0.\lb{bHtrsvs}\ee It
follows from $\bar T\neq0$ and (\rf{bHtrsvs}) that $\hat T\neq0$ and
(\rf{defts}) is well-defined. By (\rf{XFx}) with $F$ replaced by
$\bar H$, we have \be X_{\bar H}(x)=(\bar H'(x)\cdot\xi(x))X_{\hat
H}(x),\ \fa x\in D_0.\lb{XbHD0}\ee Since $\bar H(\bar
x(t)+t(0,k))\ev M$, we have \be \hat x(s)+s(0,k)=\bar
x(t(s))+t(s)(0,k)\in D_0,\ \fa s\in\R.\lb{ys+D0}\ee Then
\[\balgd\frac{\rd}{\rd s}(\hat x(s)+s(0,k))&=\frac{\rd t}{\rd s}\frac{\rd}{\rd
t}(\bar x(t)+t(0,k))
=\frac{\rd t}{\rd s}\bar TX_{\bar H}(\bar x(t)+t(0,k))\\
&=\frac{\rd t}{\rd s}\bar T\big(\bar H'(\bar x(t)+t(0,k))\cdot
\xi(\bar x(t)+t(0,k))\big)X_{\hat H}(\bar x(t)+t(0,k))\\
&=\hat TX_{\hat H}(\hat x(s)+s(0,k)).\ealgd\] The third equality
follows from (\rf{XbHD0}) and (\rf{ys+D0}). The last equality
follows from (\rf{defts}). By (\rf{defts}), the inverse function
$s(t)$ of $t(s)$ is given by \[\lt\{\balgd\frac{\rd s}{\rd
t}&=\frac{\bar T}{\hat T}\bar H'(\bar x(t)+t(0,k))\cdot \xi(\bar
x(t)+t(0,k)),\\s(0)&=0.\ealgd\rt.\] We obtain $t(1)=1$ from
\[s(1)=\frac{\bar T}{\hat T}\int^1_0\bar H'(\bar x(t)+t(0,k))\cdot
\xi(\bar x(t)+t(0,k))\rd t=1.\] Then $\hat x(1)=\bar x(1)=\bar
x(0)=\hat x(0)$ and $\hat x$ solves (\rf{fphH}).

{\bf Step 3.} Proof of (iv).

Set \[\hat\ga(t):=\hat x(t)+t(0,k).\] Since $\hat x$ solves
(\rf{fphH}) and $k\neq0$, there exists a constant $b$ such that
$\hat H(\hat\ga(t))\ev b$ and $\hat H'\neq0$ on $\hat\ga(\R)$. By
(\rf{hH'neq0}), we have $b\in(0,\frac{1}{3}\dl)$.

For (iv.1), each element $\td x\in[\hat x]$ can be written as $\td
x=\hat x+(m',m)$, where $m'\in N_0\Z^{2n-\ell},\ m\in\Z^\ell$. Since
$\hat x$ solves (\rf{fphH}), it follows from (\rf{pexthH}) that $\td
x=\hat x+(m',m)$ solves (\rf{fphH}).

For (iv.2), since $\hat H(\hat \ga(t))\ev b\in(0,\frac{1}{3}\dl)$,
by (iii), there exists $\dl_b\in(-\frac{1}{3}\dl,\frac{1}{3}\dl)$
such that (\rf{hH-1b}) holds. Then \be\hat\ga(\R)\sbs
D_{\dl_b}+N_0\Z^{2n-\ell}\ts\{0\}.\lb{hgaR}\ee By (\rf{Dsprd}) and
(\rf{Ndltactdi}), we have
\[D_{\dl_b}\sbs\tr{int}\cO_\dl\ts\R^\ell.\] By (\rf{ept}), we have
\begin{align}&(D_{\dl_b}+(m'_1,0))\cap(D_{\dl_b}+(m'_2,0))\nnb\\
\sbs&(\tr{int}\cO_\dl\ts\R^\ell+(m'_1,0))\cap(\tr{int}\cO_\dl
\ts\R^\ell+(m'_2,0))\nnb\\=\,&\es,\ \fa m'_1,m'_2\in
N_0\Z^{2n-\ell},m'_1\neq m'_2.\lb{eptt}\end{align} Since
$\hat\ga(\R)$ is connected, by (\rf{hgaR}) and (\rf{eptt}), there
exists a unique $m'_b\in N_0\Z^{2n-\ell}$ such that the orbit
$\hat\ga$ lies on $D_{\dl_b}+{(m'_b,0)}$. Then
$x(t)+t(0,k)=\hat\ga(t)-(m'_b,0)$ lies on $D_{\dl_b}$.

For (iv.3), the conclusion follows from (\rf{dtx2}), (\rf{defg}) and
elementary calculations.

For (iv.4), since $\pi\circ w$ is a non-constant $1$-periodic orbit,
then $\td\sP$ is well-defined and there exists $n(x)\in\N$ such that
(\rf{disprd}) holds.

For (iv.5), since $x(t)+t(0,k)$ lies on $D_{\dl_b}$, then $w(t)$
lies on $D_0$. Similarly to (\rf{bHtrsvs}), we have \be
H'(w(t))\cdot\xi(w(t))\neq0.\lb{Htrsvs}\ee Then $T$ determined by
(\rf{defT}) is well-defined. By $|\dl_b|<\frac{1}{3}\dl$ and
\rP{auxfg}, we have $g(\dl_b)>0$. Since $\hat T\neq0$, by
(\rf{hHtrsvs}), we have $T\neq0$. Since $x$ solves (\rf{fphH}) and
$x(t)+t(0,k)$ lies on $D_{\dl_b}$, by (\rf{dtx2}) with
$s_1=\dl_b,s_2=0,x_1(t)=x(t),x_2(t)+t(0,k)=w(t)$, we have \be
\frac{\rd}{\rd t}w(t)=e^{-\dl_b}f'(\dl_b)\hat TX_{\hat
H}(w(t))=g(\dl_b)\hat TX_{\hat H}(w(t)).\lb{dtw}\ee The second
equality follows from (\rf{defg}). Similarly to (\rf{XbHD0}), we
have \be X_H(x)=(H'(x)\cdot\xi(x))X_{\hat H}(x),\ \fa x\in
D_0.\lb{XHD0}\ee It follows from $g(\dl_b),\hat T,T\neq0$ and
(\rf{Htrsvs}) that (\rf{defths}) is well-defined. Since $z(s)=w(t)$
lies one $D_0$, we have \[\balgd\frac{\rd}{\rd s}z(s)&=\frac{\rd\hat
s}{\rd s}\frac{\rd
t}{\rd\hat s}\frac{\rd}{\rd t}w(t)\\
&=\frac{1}{T}\frac{\rd t}{\rd\hat s}g(\dl_b)\hat TX_{\hat H}(w(t))\\
&=\frac{1}{T}\frac{\rd t}{\rd\hat s}g(\dl_b)\hat
T\frac{X_H(z(s))}{H'(w(t))\cdot \xi(w(t))}\\&=X_H(z(s)).\ealgd\] The
fourth equality follows from (\rf{dtw}). The fifth equality follows
from (\rf{XHD0}). The last equality follows from (\rf{defths}). By
(\rf{defths}), the inverse function $\hat s(t)$ of $t(\hat s)$ is
given by
\[\lt\{\balgd\frac{\rd\hat s}{\rd t}&=\frac{g(\dl_b)\hat
T}{TH'(w(t))\cdot \xi(w(t))},\\ \hat s(0)&=0.\ealgd\rt.\] We obtain
$t(1)=1$ from \[\hat s(1)=\frac{g(\dl_b)\hat
T}{T}\int^1_0\frac{\rd\tau}{H'(w(\tau))\cdot\xi(w(\tau))}=1.\]
Then \[\balgd z(T)-z(0)&=w(1)-w(0)\\
&=\phi^{-\dl_b}(x(1)+(0,k))-\phi^{-\dl_b}(x(0))\\
&=\phi^{-\dl_b}(x(0)+(0,k))-\phi^{-\dl_b}(x(0))\\
&=(0,k).\ealgd\] The third equality follows from $x(1)=x(0)$. The
last equality follows from $x(0)\in D_{\dl_b}$ and (\rf{phi-sx0m})
with $s=\dl_b$. \ePf

\bR Note that the Hamiltonian function $\br H$ defined by
(\rf{defbrH1}) or (\rf{defbrH2}) also satisfies \rP{pext} except
that it is not periodic in $p$. Similar results of \cite{CZ83}
listed in Section 3 are valid for $\br H$. However, we can only
obtain $n+1$ solutions of (\rf{fphH}) with $\hat H$ replaced by $\br
H$. Then, in section 4.1, we can only obtain
$^\#\td\cP\gsl\frac{n+1}{2}$. Therefore, periodic extension of $\br
H$ in $p$ enable us to obtain more solutions of (\rf{fphH}) and
(\rf{feH}).\eR

\bP\lb{existfphH} Assume that $H$ satisfies \hlkh~with $\ell=n$.
There exists a solution of (\rf{fphH}).\eP

\bPf By the proof of Theorem 0.1 of \cite{pF92TMA}, there exists a
solution $(\bar x,\bar T)$ of (\rf{febH}). Then by (iv) of \rP{pext}
there exists a solution of (\rf{fphH}). \ePf

\setcounter{equation}{0}
\section{Saddle point reduction}

For the Hamiltonian function $\hat H:\R^{2n}\ra\R$ and the number
$\hat T$ given by \rP{pext}, define two functionals
$\cA,\cB:E=H^1(S^1,\R^{2n})\ra\R$ by
\begin{align}\cA(x)=\,&\frac{1}{2}\<Ax,x\>_{L^2}-\cB(x),\lb{defcA}\\
\cB(x)=\,&\int^1_0(\hat T\hat H(x(t)+t(0,k))+x\cdot J(0,k))\rd
t.\lb{defcB}\end{align} Critical points of $\cA$ are exactly
solutions to the equation \be Ax-\cB'(x)=0.\lb{criteq}\ee In this
section, we follow the idea of \cite{CZ83} to reduce the functional
$\cA$ to a functional $G:Z\ra\R$ defined on a finite dimensional
subspace $Z$ of $L:=L^2(\bS^1,\R^{2n})$ such that critical points of
$G$ are in one-to-one correspondence with that of $\cA$ and hence
give 1-periodic solutions of (\rf{fphH}).

Note that the spectrum of the operator $A=-J\frac{\rd}{\rd t}:E\ra
L$ is a pure point spectrum and $\sg(A)=2\pi\Z$. Let $\{E_\lm\ |\
\lm\in\R\}$ be the spectral resolution of $A$. In view of
(\rf{hH'hH''}), choose some $C\gsl0$ such that $2\pi<2\hat
TM_2+C\notin 2\pi\Z=\sg(A)$. We can define a orthogonal projection
$P$ of $L$ by \be P:=\int^{2\hat TM_2+C}_{-2\hat TM_2-C}\rd
E_\lm.\lb{defP}\ee Set $P^\bot=1-P,\ Z=P(L)$ and $Y=P^\bot(L)$. Then
$L=Z\oplus Y$ and $\dim Z=2n+2d$ for some $d\in\N$, since $2\hat T
M_2+C>2\pi$. Set $A_0=A|_Y$, $z=Px$ and $y=P^\bot x,\ x\in E$.
Equation (\rf{criteq}) is
equivalent to \begin{align} Az-P\cB'(z+y)&=0,\lb{criteq1}\\
Ay-P^\bot\cB'(z+y)&=0\ \lra\
y=A^{-1}_0P^\bot\cB'(z+y).\lb{criteq2}\end{align} The right hand
side of (\rf{criteq2}) is a contraction operator on $L$ with
contraction constant $\frac{1}{2}$. For fixed $z\in Z$, the equation
(\rf{criteq2}) has a unique solution $y=y(z)\in E$. Note that $x$
can be uniquely split as $x(z)=z+y(z)$ and equation (\rf{criteq1})
becomes $Az-P\cB'(x(z))=0$, which is equivalent to
$Ax(z)-\cB'(x(z))=0$. Define \be G:
Z\ra\R,z\mapsto\cA(x(z)).\lb{defG}\ee We have \be
G'(z)=Az-P\cB'(x(z))=Az-\hat TP\hat
H'(x+t(0,k))-J(0,k).\lb{ideG'}\ee In view of (\rf{pexthH}), from the
uniqueness of solution of (\rf{criteq2}) follows that \[
y(z+(m',m))=y(z),\ \fa z\in Z,\ (m',m)\in
N_0\Z^{2n-\ell}\ts\Z^\ell.\] Hence \begin{align}
x(z+(m',m))&=x(z)+(m',m),\tr{\quad and }\nnb\\
G'(z+(m',m))&=G'(z),\ \fa z\in Z,\ (m',m)\in
N_0\Z^{2n-\ell}\ts\Z^\ell.\lb{G'prd}\end{align} Note that if $z$ is
a critical point of $G$, the equivalent class
\[[z]:=z+N_0\Z^{2n-\ell}\ts\Z^\ell\] is a group of critical points of
$G$ with consistent energy. Moreover, the equivalent class
\[[x]:=[x(z)]:=x(z)+N_0\Z^{2n-\ell}\ts\Z^\ell\] is a group of 1-periodic
solutions of (\rf{fphH}) with consistent energy and uniquely
determined by the equivalent class $[z]$. Now we write $z\in Z$ as
$z=w+\xi$, where $w=\int^1_0z(t)\tt{d}t$. We have $w\in\ker(A)$ and
$\xi\in(\ker(A))^\bot\cap Z$. Writing $z=(w,\xi)$, we conclude from
(\rf{G'prd}) that $G'(z)=G'(w,\xi)$ is a vector field on
$T^{2n}_\ell\ts\R^{2d}$, where $T^{2n}_\ell$ denotes the
(non-standard) torus
$T^{2n}_\ell:=\R^{2n}/(N_0\Z^{2n-\ell}\ts\Z^\ell)\cong\R^{2n}/\Z^{2n}=T^{2n}$.
In summary, we have

\bP\lb{rstslt} {\bf (Lemma 1 of \cite{CZ83})} Rest points of
Lipschitz continuous vector field $G'(z)=G'(w,\xi)$ on
$T^{2n}_\ell\ts\R^{2d}$ and hence equivalent classes $[z]$s, which
consist of critical points of the $C^2$ function $G:Z\ra\R$ defined
by (\rf{defG}), are in one-to-one correspondence with equivalent
classes $[x]$ consisting of solutions of (\rf{fphH}).\eP

To study rest points of $G'$, we consider the gradient flow
\be\frac{\rd}{\rd s}z=G'(z)\lb{gradfl}\ee on
$T^{2n}_\ell\times\R^{2d}$.

\bP\lb{multrst} {\bf (Corollary of Theorem 5 of \cite{CZ83})} Denote
by $R(G)$ the set of rest points of gradient flow (\rf{gradfl}). If
$R(G)\neq\es$, then $^\#R(G)\gsl 2n+1$. \eP

\bR Propositions \rf{rstslt}-\rf{multrst} are basically the same
with corresponding results in \cite{CZ83} except that the potential
$\Phi$ in \cite{CZ83} is periodic in $2n$ directions, while in this
section the potential $\cB$ defined by (\rf{defcB}) is only periodic
in $2n-1$ directions since $k\neq0$. However, gradients of $\Phi$
and $\cB$ are all periodic in $2n$ directions. Hence, the idea of
\cite{CZ83} is valid for the multiplicity proof of solutions of
(\rf{fphH}).\eR

\bR\lb{2n+1rst} By Propositions \rf{existfphH} and
\rf{rstslt}-\rf{multrst}, if $H$ satisfies \hlkh~with $\ell=n$,
there exist $(x_i,b_i,\dl_{b_i}),\ 1\lsl i\lsl 2n+1$ such that $x_i$
solves (\rf{fphH}), $b_i\in(0,\frac{1}{3}\dl),\
\dl_{b_i}\in(-\frac{1}{3}\dl,\frac{1}{3}\dl),\ x_i(t)+t(0,k)$ lies
on $D_{\dl_{b_i}}$ and $x_i$s are geometrically distinct in the
sense of (\rf{gd}). In the next section, these $2n+1$ solutions are
classified into two categories and give at least $n+1$ geometrically
distinct solutions of (\rf{feH}).\eR

\setcounter{equation}{0}
\section{Proof of Theorems}

\subsection{Proof of Theorem \rf{rotfe}}

\bP\lb{pfgd} Let $k_i\in\Z^\ell\bs\{0\}$ for $i\in\{1,2\}$. Assume
that $(x_i,\tau_i)$ solves (\rf{fphH}) with $(k,\hat T)$ being
replaced by $(k_i,\tau_i)$ and $x_i(s)+s(0,k_i)$ lies on $D_{\dl_i}$
with $\dl_i\in(-\frac{1}{3}\dl,\frac{1}{3}\dl)$. If
\[w_i(s):=\phi^{-\dl_i}(x_i(s)+s(0,k_i)),\ i=1,2,\] are geometrically
the same in the sense of \be\pi\circ w_1(\R)=\pi\circ
w_2(\R),\lb{geosam}\ee where $\pi$ denotes the covering projection
$\R^{2n}\ra\R^{2n-\ell}\ts\R^\ell/\Z^\ell$. Then

(i) There exists a reparametrization $s:\R\ra\R,\ t\mt s(t)$ such
that \begin{align}\dot
s(t)&\ev\frac{g(\dl_1)\tau_1}{g(\dl_2)\tau_2},\ \fa t\in\R,\lb{idedts}\\
\pi\circ w_1(t)&=\pi\circ w_2(s(t)),\ \fa t\in\R,\lb{piw1=piw2st}\\
k_1&=\dot s(0)k_2,\lb{idedts0}\end{align} where $g:\R\ra\R$ is
defined by (\rf{defg}).

(ii) There exists $\ep\in\{\pm1\}$ such that
\be\frac{g(\dl_1)\tau_1}{g(\dl_2)\tau_2}=
\ep\frac{|k_1|}{|k_2|}.\lb{idegk}\ee\eP

\bPf For (i), without loss of generality, we assume that \be\pi\circ
w_1(0)=\pi\circ w_2(0).\lb{piw10=piw20}\ee Since $(x_i,\tau_i)$
solves (\rf{fphH}) with $(k,\hat T)$ being replaced by
$(k_i,\tau_i)$ and $x_i(s)+s(0,k_i)$ lies on $D_{\dl_i}$ with
$\dl_i\in(-\frac{1}{3}\dl,\frac{1}{3}\dl)$, it follows from (v.3) of
\rP{pext} that $\pi\circ w_i$ solves \be\lt\{\balgd\frac{\rd}{\rd
t}\pi\circ w_i(t)&=g(\dl_i)\tau_iX_{\hat H_\pi}(\pi\circ w_i(t)),\\
\pi\circ w_i(1)&=\pi\circ w_i(0).\ealgd\rt.\lb{piwisol}\ee Define
$s:\R\ra\R$ by \be s(t):=\dot
s(0)t:=\frac{g(\dl_1)\tau_1}{g(\dl_2)\tau_2}t,\ \fa
t\in\R.\lb{defstt}\ee Set \be\bar w_1(t):=w_2(s(t)).\lb{defbw1}\ee
Then \begin{align}\frac{\rd}{\rd t}\pi\circ\bar w_1(t)=\,&\frac{\rd
s}{\rd t}\frac{\rd}{\rd s}\pi\circ w_2(s)\nnb\\=\,&\dot
s(t)g(\dl_2)\tau_2X_{\hat H_\pi}(\pi\circ
w_2(s))\nnb\\=\,&g(\dl_1)\tau_1X_{\hat H_\pi}(\pi\circ\bar w_1(t)),\lb{pibw1sol1}\\
\pi\circ\bar w_1(0)=\,&\pi\circ w_2(0)=\pi\circ
w_1(0).\lb{pibw1sol2}\end{align} The second equality follows from
(\rf{piwisol}). The third equality follows from (\rf{defstt}). By
(\rf{piwisol}) with $i=1$, (\rf{pibw1sol1}) and (\rf{pibw1sol2}),
the basic uniqueness theorem for the initial value problem of
ordinary differential equations yields \be\pi\circ
w_1(t)=\pi\circ\bar w_1(t),\ \fa t\in\R.\lb{piw1=pibw1}\ee Then
(\rf{piw1=piw2st}) follows from (\rf{defbw1}) and (\rf{piw1=pibw1}).

By (\rf{idedts}) and (\rf{piw1=piw2st}), we have \[\balgd\pi\circ
w_2(s+\dot s(0))&=\pi\circ w_2(s(t)+\dot s(0))=\pi\circ w_2(\dot
s(0)(t+1))=\pi\circ w_2(s(t+1))\\
&=\pi\circ w_1(t+1)=\pi\circ w_1(t)=\pi\circ w_2(s(t))=\pi\circ
w_2(s).\ealgd\] This proves that $\dot s(0)$ is a period of
$\pi\circ w_2$. Since $(x_2,\tau_2)$ solves (\rf{fphH}) with
$(k,\hat T)$ being replaced by $(k_2,\tau_2)$ and $x_2(s)+s(0,k_2)$
lies on $D_{\dl_2}$ with $\dl_2\in(-\frac{1}{3}\dl,\frac{1}{3}\dl)$,
by (v.5) of \rP{pext}, there exists $m_2\in\Z\bs\{0\},\ n(x_2)\in\N$
such that \be\dot s(0)=\frac{m_2}{n(x_2)}.\lb{idedts0-1}\ee By
(\rf{idedts}), (\rf{piw1=piw2st}) and (\rf{idedts0-1}), there exists
$j\in\Z^\ell$ such that
\be\phi^{-\dl_1}(x_1(t)+t(0,k_1))=\phi^{-\dl_2}\lt(x_2\lt
(\frac{m_2}{n(x_2)}t\rt)+\frac{m_2}{n(x_2)}t(0,k_2)\rt)+(0,j).
\lb{phi-dl1=phi-dl2}\ee Taking the difference of
(\rf{phi-dl1=phi-dl2}) with $t=0$ and $n(x_2)$, we have
\begin{align} n(x_2)(0,k_1)&=\phi^{-\dl_1}(x_1(0)+n(x_2)(0,k_1))-
\phi^{-\dl_1}(x_1(0))\nnb\\
&=\phi^{-\dl_1}(x_1(n(x_2))+n(x_2)(0,k_1))-
\phi^{-\dl_1}(x_1(0))\nnb\\
&=\phi^{-\dl_1}(x_1(t)+t(0,k_1))\bigg|^{t=n(x_2)}_{t=0}\nnb\\
&=\phi^{-\dl_2}\lt(x_2\lt
(\frac{m_2}{n(x_2)}t\rt)+\frac{m_2}{n(x_2)}t(0,k_2)\rt)\bigg|^{t=n(x_2)}_{t=0}
\nnb\\&=\phi^{-\dl_2}(x_2(m_2)+m_2(0,k_2))-\phi^{-\dl_2}(x_2(0))\nnb\\
&=\phi^{-\dl_2}(x_2(0)+m_2(0,k_2))-\phi^{-\dl_2}(x_2(0))\nnb\\
&=m_2(0,k_2).\lb{nx20k1=}\end{align} The first and the last
equalities follow from $x_i(0)\in D_{\dl_i},\
n(x_2)k_1,m_2k_2\in\Z^\ell$ and (\rf{phi-sx0m}) with $s_i=\dl_i$
respectively. The second and the fifth equalities follow from
$n(x_2)\in\N,\ m_2\in\Z\bs\{0\}$ and $x_i$ is $1$-periodic. Then
(\rf{idedts0}) follows from (\rf{idedts0-1}) and (\rf{nx20k1=}).

For (ii), the conclusion follows from (\rf{idedts}) and
(\rf{idedts0}). \ePf

Given $k\in\Z^n\bs\{0\}$, by \rR{2n+1rst}, we obtain at least $2n+1$
solutions of (\rf{fphH}), denoted by $x_1,x_2,\cd,x_{2n+1}$, which
are distinct in the sense of (\rf{gd}). There exist
$\dl_i\in(-\frac{1}{3}\dl,\frac{1}{3}\dl),1\lsl i\lsl2n+1$ such that
\be w_i(t):=\phi^{-\dl_i}(x_i(t)+t(0,k))\in D_0,\ \fa t\in\R,\ 1\lsl
i\lsl 2n+1.\lb{defwit}\ee Set \begin{align}\cP_1&:=\{w_i\ |\
\dl_i\in(-\frac{1}{3}\dl,-\frac{1}{3}\dl+\dl^+],\
1\lsl i\lsl 2n+1\},\lb{defcP1}\\
\cP_2&:=\{w_i\ |\ \dl_i\in[-\frac{1}{3}\dl+\dl^+,\frac{1}{3}\dl),\
1\lsl i\lsl 2n+1\},\lb{defcP2}\end{align} where $\dl^+$ is defined
by (\rf{defdlta+}). Recall that $g$ defined by (\rf{defg}) is
positive on $(-\frac{1}{3}\dl,\frac{1}{3}\dl)$, increases strictly
on $(-\frac{1}{3}\dl,-\frac{1}{3}\dl+\dl^+]$ and decreases strictly
on $[-\frac{1}{3}\dl+\dl^+,\frac{1}{3}\dl)$. We have

\bP\lb{pfgs} For every $i\in\{1,2\}$, elements of $\cP_i$ are
pair-wise distinct in the sense of (\rf{gd}). \eP

\bPf We argue by contradiction. Suppose that there exists a pair
$w_{i_1},\ w_{i_2}$ belonging to $\cP_i$ for some $i\in\{1,2\}$ such
that $\pi\circ w_{i_1}(\R)=\pi\circ w_{i_2}(\R)$. By (iii) of
\rP{pfgd} with $(k_1,\tau_1)=(k_2,\tau_2)=(k,\hat T)$, we have
$g(\dl_{i_1})=\ep_1g(\dl_{i_2})$. By \rP{auxfg}, we have
$g(\dl_{i_1})=g(\dl_{i_2})>0$ and hence $\dl_{i_1}=\dl_{i_2}$. By
(\rf{idedts}), we have $\dot s(0)=1$. By (\rf{idedts0-1}) and
(\rf{phi-dl1=phi-dl2}), we have \begin{align}
x_{i_1}(t)+t(0,k)&=\phi^{\dl_{i_1}}\big(\phi^{-\dl_{i_2}}
(x_{i_2}(t)+t(0,k))+(0,j)\big)\nnb\\
&=\phi^{\dl_{i_1}}\circ\phi^{-\dl_{i_2}}(x_{i_2}(t)+t(0,k))+(0,j)\nnb\\
&=x_{i_2}(t)+t(0,k)+(0,j).\lb{xi1=xi2}\end{align} The second
equality follows from (\rf{defwit}) and (\rf{phi-sx0m}) with
$s=\dl_{i_1}$. The last equality follows from $\dl_1=\dl_2$.
However, the equality (\rf{xi1=xi2}) contradicts the assumption
$\pi\circ x_{i_1}(\R)\neq\pi\circ x_{i_2}(\R)$. \ePf

{\bf Proof of Theorem 1.1.} By \rP{pfgs}, we have
\[\balgd^\#\{w_i,\ 1\lsl i\lsl 2n+1\ &|\ w_i\tr{s are geometrically
distinct in the sense of (\rf{gd})}.\}\\
&\gsl\max\{^\#\cP_1,^\#\cP_2\}\gsl\frac{2n+1}{2}>n.\ealgd\] By
(iv.5) of \rP{pext}, we have $\td\cP(M,k)\gsl n+1$. \QED

\subsection{Proof of \rT{alter}}

Since $H^{-1}(M)$ is contact, by (\rf{defbrH1}) or (\rf{defbrH2}),
we obtain a Hamiltonian function $\br
H:\R^{2n-\ell}\ts(\R^\ell/\Z^\ell)\ra\R$ which is bounded in
$p\in\R^{2n-\ell}$ and possesses $H^{-1}(M)$ as a regular energy
hypersurface. Since $2n-\ell\gsl 1$, by Theorem 1.2 of \cite{mJ94}
and (v.3) of \rP{pext}, there exists a nonconstant solution $(z,T)$
of (\rf{feH}) with some $k_0\in\Z^\ell$. If $k_0=0$ then
$^\#\td\cP(M,0)\gsl1$. If $k_0\neq0$, we extend $\br H$ periodically
in $p$ to obtain a new Hamiltonian function $\hat H$, which
satisfies \rP{pext}. Repeating the proof of
Section 4.1, we have $^\#\td\cP(M,k_0)\gsl n+1$. \QED\\

ACKNOWLEDGEMENTS.\quad This article is a part of my Ph. D. thesis. I
would like to thank my thesis advisor, Professor Yiming Long, for
introducing me to this filed, drawing my attention to the work of P.
Felmer, and many helpful discussions on this topic. I would also
thank the referee for his/her valuable comments on the paper.

\clearpage

\bfg\cpts{1}\centering
\scalebox{0.5}{\rotatebox{-90}{\includegraphics{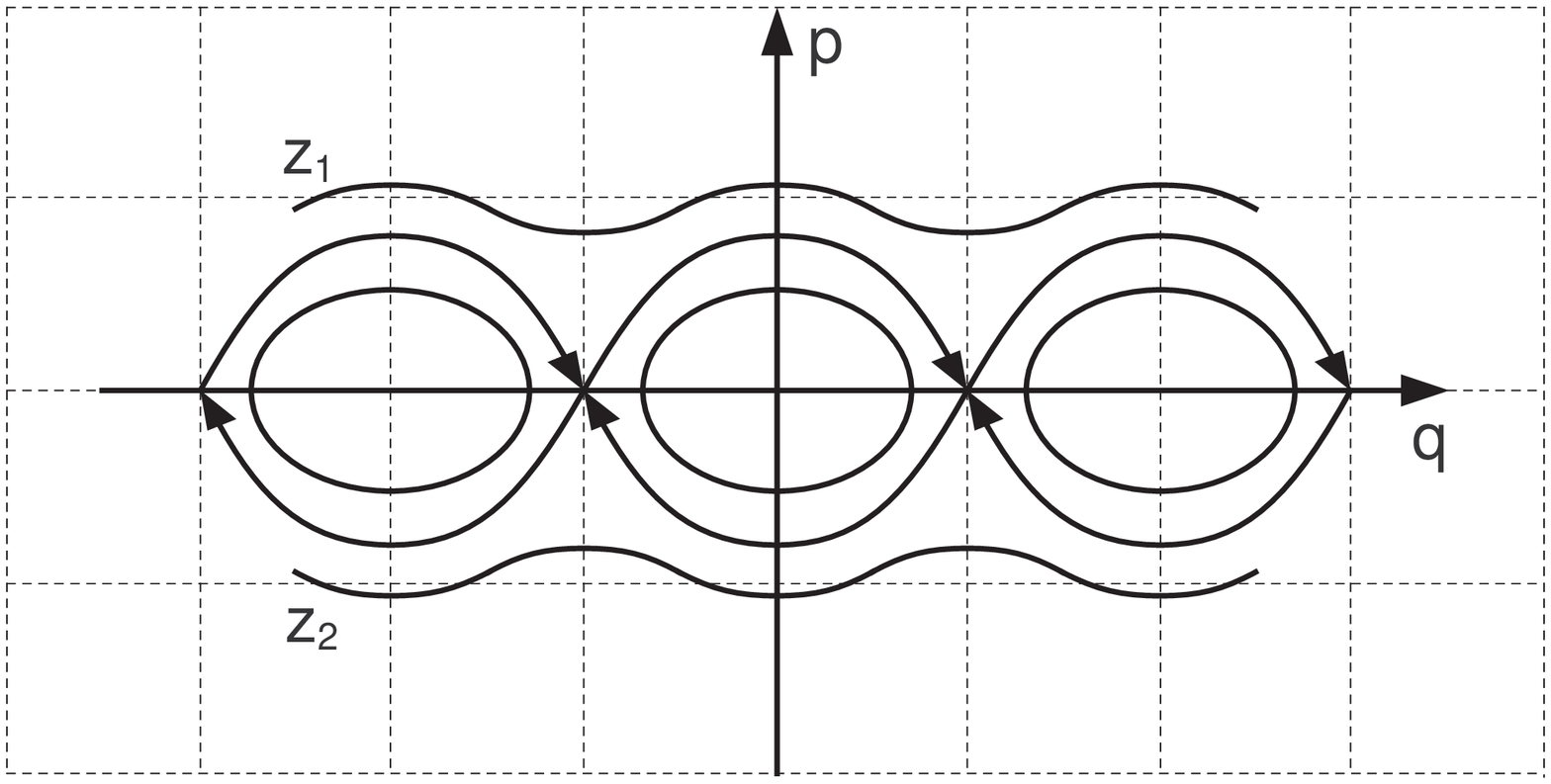}}}
\caption{simple pendulum}\lb{fig.1}\efg

\clearpage

\bfg\cpts{1}\centering
\hs{-3.5cm}\scalebox{0.5}{\rotatebox{-90}{\includegraphics{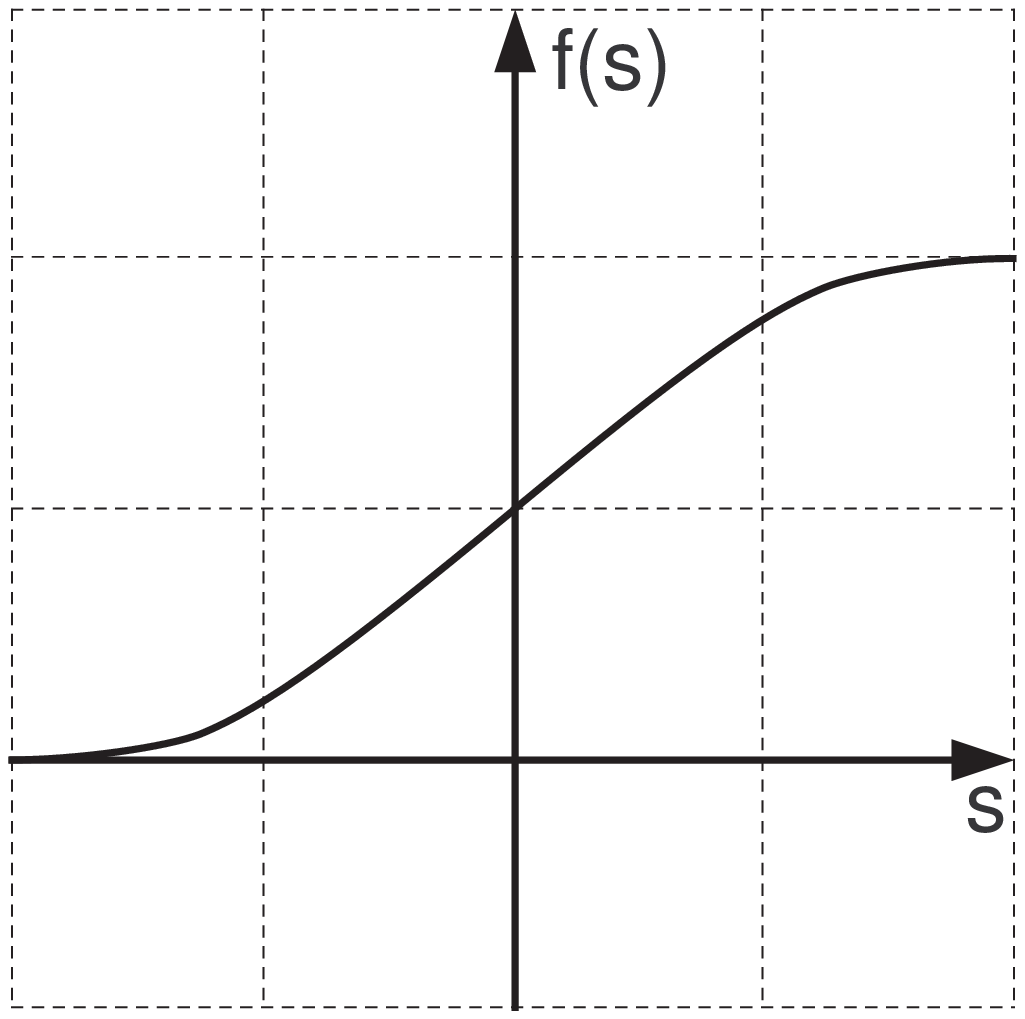}}}
\caption{auxiliary function $f$}\lb{fig.f}\efg

\clearpage

\bfg\cpts{1}\centering
\scalebox{0.5}{\rotatebox{-90}{\includegraphics{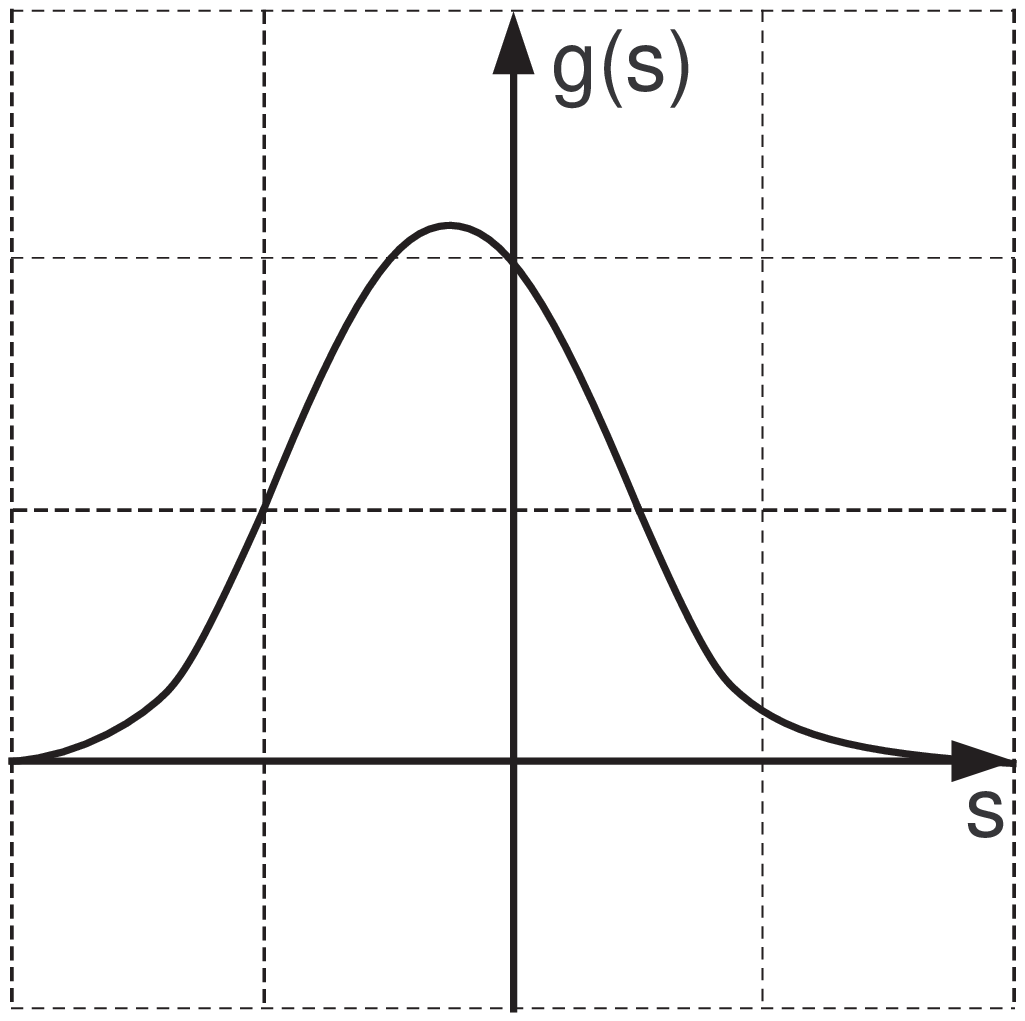}}}
\caption{auxiliary function $g$}\lb{fig.g}\efg

\end{document}